\documentclass[12pt]{article}

\usepackage{amsmath,amssymb}
\usepackage{pstricks,pst-node,graphicx}

\newtheorem{theorem}{Theorem}[section]

\newtheorem{example}{Example}[section]

\newcommand{\be}{\begin{equation}}
\newcommand{\ee}{\end{equation}}
\newcommand{\eea}{\end{eqnarray}}
\newcommand{\bea}{\begin{eqnarray}}

\begin{document}
\title{{\bf \Large An ergodic diffusion with unbounded inward and outward drift}}

\author{{\bf Horst Thaler} \\[1ex] {\small Department of Mathematics and Informatics,
University of Camerino,} \\
{\small Via Madonna delle Carceri 9,
I--62032, Camerino (MC), Italy;}\\
{\small
E-mail:
horst.thaler@unicam.it}}

\date{}

\maketitle

{\abstract It is argued that a diffusion may be ergodic even though the drift field has unbounded outward-directed parts. The discussion employs stochastic and numerical methods.
\\[1ex]
}
\section{Introduction}
We shall investigate the ergodicity of the following stochastic differential equation (SDE)
\begin{equation}\label{IO}
\left\{
\begin{aligned}
dX_t &=b(X_t)dt+dW_t, \\
X_0 &=\zeta, &
\end{aligned}
\right.
\end{equation}
where $\zeta\in \mathbb{R}^2$ and $(W_t)_{t\geq 0}$ is a $2$-dimensional Brownian motion.
The drift $b:\mathbb{R}^2\rightarrow \mathbb{R}^2$ is given by
\be\label{driftz4}
b_1=-4x_1^3+12x_1x_2^2, \quad b_2=-12x_1^2x_2+4x_2^3.\\[1ex]
\ee
\begin{figure}[h] \label{fig1}
\begin{center}
\includegraphics[scale = 0.7]{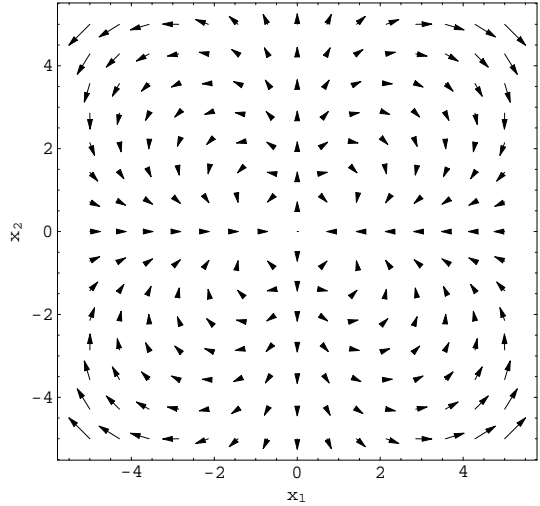}
\caption{\label{driftz4gr} The drift field $b$ according to (\ref{driftz4})}
\end{center}
\end{figure}
Here a stochastic process $X=(X_t)_{t\geq 0}$ is termed ergodic if there exists a unique invariant probability measure $\mu$ for $X$ and if for every $f\in L^1(\mathbb{R}^d,\mu)$ it holds that
\be\label{ergf}
\lim_{T\rightarrow \infty}\frac{1}{T}\int_0^T f(X_t)dt = \int_{\mathbb{R}^d}f d\mu, \quad P_\zeta-\mathrm{a.s.}, \;\mathrm{for\;all}\; \zeta \in \mathbb{R}^d,
\ee
where $P_\zeta$ denotes the law of $X$ with starting point $\zeta$.
Condition (\ref{ergf}) is also called the law of large numbers of the process $X$. Under very general assumptions it can be shown that the existence of an invariant probability measure is equivalent to $X$ being positive recurrent, see \cite[Theorem 4.1]{Kli}.
As we shall see it is difficult to classify the diffusion $X$, which solves equation (\ref{IO}), as being transient, recurrent or positive recurrent. This is due to the funny drift $b$ itself.
Before proceeding we fix some notation.

{\it Notation:} $\mathbb{E}_\zeta[\,\cdot\,]$ shall denote the expectation with respect to the law $P_\zeta$. The standard inner product of two vectors $x,x'\in \mathbb{R}^d$ is written as $x\cdot x'$. We set $r=r(x):=\sqrt{x\cdot x}$. Furthermore, $B(\mathbb{R}^d)$ is the space of bounded, real-valued measurable functions and $C^2(\mathbb{R}^d)$ stands for the space of twice continuously differentiable functions. The indicator function of a set $A$ will be denoted by $1_A$.

Let us have a closer look at $b$ defined in (\ref{driftz4}). Let $e_r:=(x_1,x_2)/r$ be the radial unit vector field, then we may distinguish two cases
\be
\begin{aligned}
(a)\;\; e_r\cdot b &< 0, \;\mathrm{if}\;\; |x_2|<|x_1|, \\
(b)\;\; e_r\cdot b &> 0, \;\mathrm{if}\;\; |x_2|>|x_1|.
\end{aligned}
\ee
Therefore, in case ($a$) we have an inward drift and in case ($b$) an outward drift. Moreover, evaluating $e_r\cdot b$ with $x_1= r\cos\varphi,\, x_2=r\sin\varphi,\,\varphi\in [-\pi,\pi],$  we see that in both cases there results a growth behavior proportional to $r^3$ and that the proportionality constant is positive and increases as $|\varphi|\rightarrow \pi/2$ for $\{\pi/4<|\varphi|<3\pi/4\}$ and decreases, being negative, as $|\varphi|\rightarrow 0$ or $|\varphi|\rightarrow \pi$ for $\{|\varphi| <\pi/4\}$ and $\{|\varphi| > 3\pi/4\}$, respectively. It equals $0$ for $|\varphi|=\pi/2,|\varphi|=3\pi/4.$ This means that moving away from the origin there are areas with unbounded inward drift as well as areas with unbounded outward drift. Moreover, it should be noted that $\mathrm{curl}(b_1,b_1)=-12x_1x_2$ which shows that it is impossible to write the drift as $b=-\nabla V$ for some differentiable function $V$.

The paper is organized as follows. In section \ref{trre} we collect some background material on the different large time asymptotics that a stochastic process may show. These properties are discussed in terms of the coefficients of the generator of $X$. In section \ref{diundr} we return to the investigation of the diffusion solving (\ref{IO}). It is made evident that the latter diffusion cannot be handled using the standard criteria of section \ref{trre}. Finally, the numerical results for the long-term behavior of (\ref{ergf}) are presented which give strong indications that the diffusion solving (\ref{IO}) is ergodic.

\section{Criteria for transience and recurrence}\label{trre}

Here we discuss criteria that allow for a complete classification of a diffusion, at least if the latter are applicable. To this end we consider a class of more general SDEs, namely
\begin{equation}\label{lequ}
\left\{
\begin{aligned}
dX_t &= b(X_t)dt+dW_t, \\
X_0 &= \zeta,
\end{aligned}
\right.
\end{equation}
where $(W_t)_{t\geq 0}$ is a $d$-dimensional Brownian motion and $b:\mathbb{R}^d \rightarrow \mathbb{R}^d$ is a smooth vector field. The semigroup $(S_t)_{t\geq 0}, \; S_t:B(\mathbb{R}^d)\rightarrow B(\mathbb{R}^d)$, associated with $X$ is given by
$S_tf(\zeta):=\mathbb{E}_\zeta\left[f(X_t)1_{\{t<\eta\}}\right]$, where $\eta$ denotes the explosion time of $X$.
We shall assume that the generator $L$ of $X$ uniquely determines $(S_t)_{t\geq 0}$ and that on $C^2(\mathbb{R}^d)\subset D(L)$ it is given by
$$L=\textstyle{\frac{1}{2}}\Delta+ b\nabla,\quad \mathrm{with\;formal\;adjoint\;\;}L^\ast(\cdot)=
\mathrm{div}((\textstyle{\frac{1}{2}}\nabla-b)(\cdot)).$$
Since a priori we do not demand the semigroup to be strongly continuous the generator $L$ might be given also in the weak sense, see e.g. \cite[Section 2.2]{Dyn} and \cite[Section 2.3]{Lor}. From now on processes sharing the foregoing properties will simply be called diffusions.

In \cite{Bha} Bhattacharya has given quite general sufficient criteria in terms of the coefficients $b_i$ that allow to decide whether a diffusion is transient, recurrent or that it even admits a finite invariant measure. We recall that the diffusion $X$ is said to be recurrent if for any $\zeta\in \mathbb{R}^d$ and any open set $U\subset \mathbb{R}^d$ we have
$$P_\zeta(X_{t_n}\in U \; for \; a \; sequence \;of \; {\rm {\it t_n}'}s \; increasing \; to \; \infty)=1.$$
If for any $\zeta\in \mathbb{R}^d$ it holds that
$$P_\zeta(\lim_{t\rightarrow\infty}r(X_t)=\infty)=1,$$
then $X$ is called a transient process. Finally, $X$ is said to be positive recurrent, if $\mathbb{E}_\zeta[\sigma_U]<\infty$, for all $\zeta\in \mathbb{R}^d$ and all open sets $U\subset \mathbb{R}^d$, where $\sigma_U$ is the first entrance time of $U$.

Adopting the same notation as in \cite{Bha}, we introduce the following functions
\begin{equation}
\begin{aligned}
C(x)\, &=\, 2 \sum_{i=1}^d x_i b_i(x) \\
\overline{\beta}(r)\, &=\, \sup_{|x|=r}(d-1+C(x)),\quad \underline{\beta}(r) = \inf_{|x|=r}(d-1+C(x))\\
\overline{I}(r) \, &=\, \int_{r_0}^r \frac{\overline{\beta}(u)}{u} du, \quad \underline{I}(r) = \int_{r_0}^r \frac{\underline{\beta}(u)}{u} du
\end{aligned}
\end{equation}

With these preparations we may quote the following two theorems from \cite{Bha}

\begin{theorem}

(a) If for some $r_0>0$
\begin{equation}\label{cr1}
\int_{r_0}^\infty \exp\{- \overline{I}(u)\}du =\infty,
\end{equation}
then the diffusion with generator $L$ is recurrent. \\
(b) If for some $r_0>0$
\begin{equation}\label{cr2}
\int_{r_0}^\infty \exp\{- \underline{I}(u)\}du <\infty,
\end{equation}
then the diffusion with generator $L$ is transient.

\end{theorem}

\begin{theorem}

(a) The diffusion with generator $L$ is recurrent and admits a finite invariant measure (unique up to a constant multiple) if there exists a $r_0>0$ such that (\ref{cr1}) holds and
\be\label{cr4}
\int_{r_0}^\infty \exp\{ \overline{I}(u)\}du <\infty.
\ee
(b) If there exists a $r_0>0$ such that (\ref{cr1}) holds and
\be\label{cr5}
\lim_{N\rightarrow \infty} \left.\int_{r_0}^N \exp\{- \overline{I}(s)\}\left(\int_{r_0}^s \exp\{ \overline{I}(u)\}du\right)ds \right\slash \int_{r_0}^N \exp\{- \underline{I}(u)\}du =\infty,
\ee
then the recurrent diffusion does not admit a finite invariant measure.

\end{theorem}

Henceforth a unique finite invariant measure will always be normalized to become a probability measure.
In the following examples we explore some possible scenarios. For simplicity only radially symmetric problems are considered and we set $I(r):=\overline{I}(r)=\underline{I}(r)$.

\begin{example}

(a) Let $b=0$ in (\ref{lequ}), then $X_t=W_t$ and an elementary calculation gives by (\ref{cr1}) and (\ref{cr2}) that Brownian motion is recurrent in $d=1,2$ and transient for $d\geq 3.$ \\[1ex]
(b) Let $V(x)=-r^{-\alpha}$ and $b= -\nabla V$, with $\alpha>0$. In this case we find
$$
I(r)=\int_{r_0}^r\frac{(d-1-2\alpha u^{-\alpha})}{u}du=(d-1)(\ln r-\ln r_0)+2(r^{-\alpha}-r_0^{-\alpha}).
$$
By elementary calculations one deduces that (\ref{cr1}) holds for $d=1,2$, implying recurrence, and that (\ref{cr2}) is valid for $d\geq 3$ entailing transience.
On the other hand, evaluating the quotient of the l.h.s in (\ref{cr5}) gives the following estimate
$$
{quotient\; of\; integrals \; of}\; {\it (\ref{cr5})} \geq C(r_0)N,
$$
for $N$ large enough, with $C(r_0)$ some positive constant. This
shows that there is no finite invariant measure for any dimension. Still there is a $\sigma$-finite invariant measure $\nu(dx)=e^{-2V}dx$ for every dimension, since the stationary Fokker-Planck equation has a solution
$$
L^\ast(e^{-2V})=
\mathrm{div}((\textstyle{\frac{1}{2}}\nabla-b)e^{-2V})=0.
$$
(c) Let $V(x)=r^\alpha$ and $b=-\nabla V$, where $\alpha >0.$ Here we obtain
$$I(r)=\int_{r_0}^r \left(\frac{(d-1)}{u}-2\alpha u^{\alpha-1}\right)du=\ln{r^{d-1}}-2 r^\alpha-\ln{r_0^{d-1}}+2 r_0^\alpha,$$
which entails properties (\ref{cr1}) and (\ref{cr4}). Now the diffusion is recurrent with finite invariant measure.
\end{example}

{\it Remark.} The functions $V$ in examples ($b$) and ($c$) are not smooth at the origin. But the latter may be easily rendered such.\\[1ex]
\indent
It is clear that the diffusive part $dW_t$ tends to spread a ``particle'' moving according to (\ref{lequ}). This spreading may be compensated by a sufficiently strong drift $b$. As can be seen from the previous example (c), if the driving field $b$ is inward directed and the gradient of a potential that grows like a positive power, then we are guaranteed an invariant probability measure $\mu$ on $\mathbb{R}^d$. In fact in the case of example (c) the diffusion is even ergodic, meaning that for every $f\in L^1(\mathbb{R}^d,\mu)$ we have
\be\label{erg}
\lim_{T\rightarrow \infty}\frac{1}{T}\int_0^T f(X_t)dt=\int_{\mathbb{R}^d}f d\mu, \quad P_\zeta-\rm{a.s}, \;\mathrm{for \; all\;} \zeta\in \mathbb{R}^d.
\ee
Notice that, for $f$ bounded and measurable, property (\ref{erg}) immediately implies
\be\label{erg1}
\lim_{T\rightarrow \infty}\frac{1}{T}\int_0^T S_tf(\zeta)dt=\int_{\mathbb{R}^d}f d\mu, \quad \; \zeta\in \mathbb{R}^d.
\ee
A large class of processes for which properties
(\ref{erg}) and (\ref{erg1}) hold is given by recurrent diffusions with invariant probability
measure that satisfy in addition the following conditions (cf. \cite[Section 1.3]{Kun})

{\it Condition} $(A)$ The semigroup $(S_t)_{t\geq 0}$ maps $C_0(\mathbb{R}^d)$ into $C_0(\mathbb{R}^d)$ and is strongly continuous.

{\it Condition} $(B)$ There is a Borel measure $\nu$ on $\mathbb{R}^d$ with support on $\mathbb{R}^d$ and a strictly positive function $p_t(x,y)$ continuous for $(t,x,y)\in(0,\infty)\times\mathbb{R}^d\times\mathbb{R}^d$ such that the transition probability can be written as $K_t(x,dy)=p_t(x,y)\nu(dy).$ \\[1ex]
\indent

\section{An ergodic diffusion with strong outward drift}\label{diundr}

We now come to the discussion of the diffusion $(X_t)_{t\geq 0}$ which solves (\ref{IO}).
The peculiar form of (\ref{driftz4}) is related to a complex differentiation, namely if $z=x_1+ix_2$ and $g=(x_1+i x_2)^4$, then
$$b_1=-\mathrm{Re}\left(\frac{d}{dz}(g)\right), \quad b_2=-\mathrm{Im}\left(\frac{d}{dz}(g)\right).$$
We know at least that $X$ has continuous sample paths, maps $B(\mathbb{R}^d)$ into $C(\mathbb{R}^d)$ and is conservative, i.e. its explosion time is almost surely infinite. The first property is a consequence of our assumptions, see e.g. \cite{Ike}. The second and third properties are consequences of Lemma 2.5 and condition (4.1) of \cite{Bha}.

Let us turn to the transience and recurrence properties of our diffusion.
The absolute maxima and minima of $(1+C(x))$ given the constraint $|x|=r$ are found to be
$\overline{\beta}(r)=1+8r^4$ and $\underline{\beta}(r)=1-8r^4$, respectively. This gives $\overline{I}(r)=\ln{r}-\ln{r_0}+ 2(r^4-r_0^4)$ and $\underline{I}(r)=\ln{r}-\ln{r_0}- 2(r^4-r_0^4)$, hence
\be
\int_{r_0}^\infty \exp\{-\overline{I}(u)\}du = r_0 e^{2r_0^4}\int_{r_0}^\infty \frac{1}{u} e^{-2u^4}du <\infty
\ee
\be
\int_{r_0}^\infty \exp\{-\underline{I}(u)\}du = r_0 e^{-2r_0^4}\int_{r_0}^\infty \frac{1}{u} e^{2u^4}du =\infty .
\ee
This shows that none of the criteria (\ref{cr1}) and (\ref{cr2}) applies. As far as the existence of a finite invariant measure is concerned we use the estimate
\be\label{est1}
\mathrm{quotient\; of\;integrals\;of} \; (\ref{cr5}) \left. \leq C(r_0)\int_{r_0}^N (s-r_0) ds \right\slash \int_{r_0}^N \frac{1}{u} e^{2u^4}du,
\ee
for some positive constant $C(r_0).$
Since the r.h.s of (\ref{est1}) tends to 0 as $N\rightarrow \infty$ criterion (\ref{cr5}) is not valid. Hence, even if the diffusion were recurrent we could not tell whether it possesses an invariant probability measure.

\begin{figure}
\begin{center}
\includegraphics[scale = 0.9]{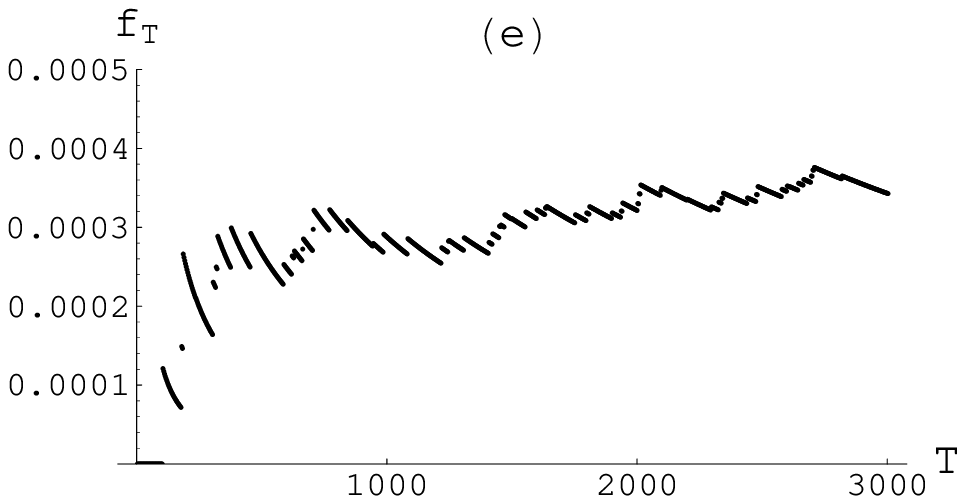}
\includegraphics[scale = 0.9]{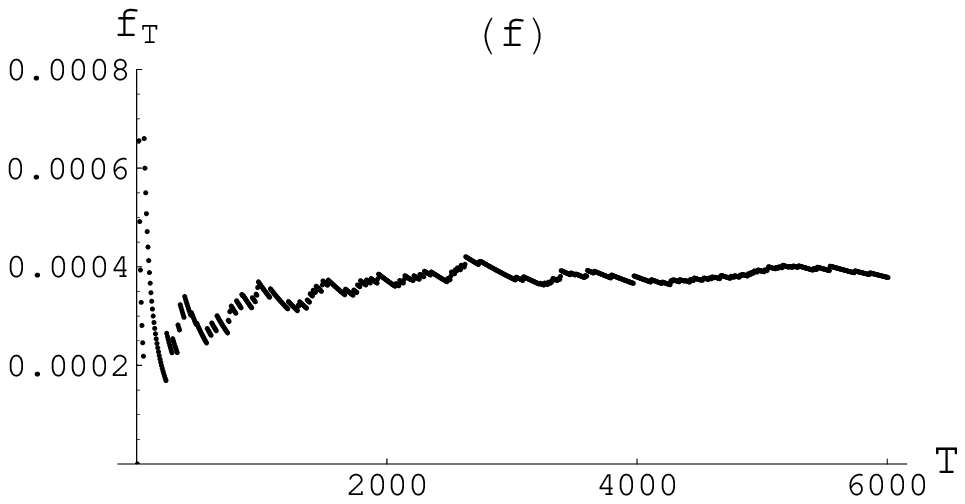}\\[2ex]
\includegraphics[scale = 0.9]{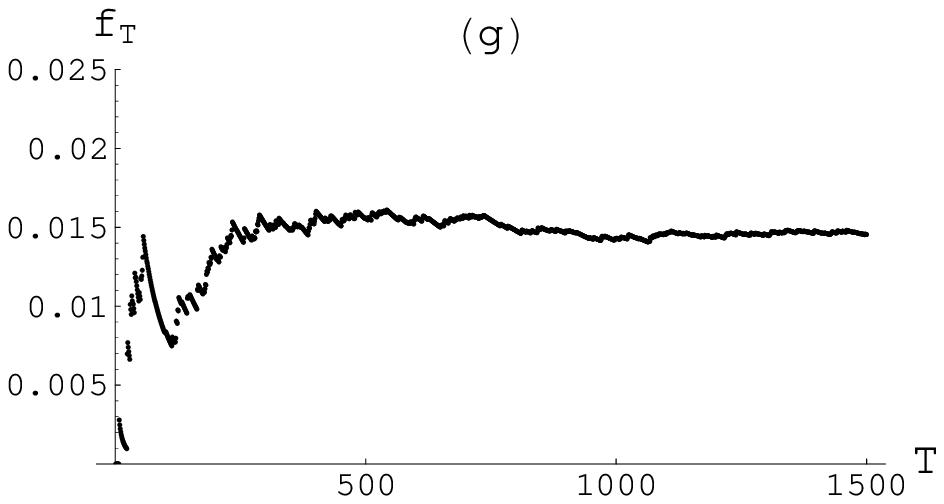}
\includegraphics[scale = 0.9]{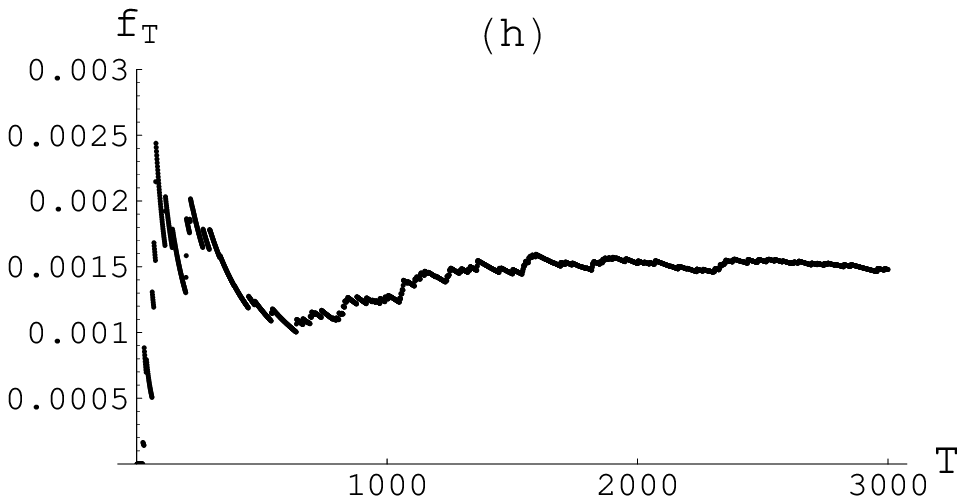}
\caption{\label{graphs} Behavior of $f_T$, with $f=1_{B_1}(x-x_0)$ where $B_1=\{x\in \mathbb{R}^2 \;|\;r(x)<1\}$ and for $x_0$ the following values have been chosen. (a) $(0,0)$, (b) $(2,0)$, (c) $(0,2)$, (d) $(3,0)$, (e) and (f) $(0,3)$, (g) 2$(1/\sqrt{2},1/\sqrt{2})$, (h) 3$(1/\sqrt{2},1/\sqrt{2})$. The starting points $\zeta$ have been picked randomly from $[-10,10]^2$.}
\end{center}
\end{figure}

In order to get some insight into the problem of the diffusion (\ref{IO}) several numerical tests have been performed whose aim was to study the behavior of
\be\label{fT}
f_T:= \frac{1}{T}\int_0^T f(X_t) dt ,
\ee
as $T$ increases.
As $f_T$ is an expression depending on the sample paths a strong Taylor approximation (of order $1.5$) was used for the latter. More precisely, the strong Taylor scheme is implemented through (cf. \cite[p. 355, eq. 4.12]{Klo})
\be\label{alg}
Y^k_{n+1}=Y_n^k + b_k \Delta + \Delta W^k + {\textstyle \frac{1}{2}} L^0 b_k \Delta^2 +  L^k b_k \Delta Z^k,\quad k=1,2,
\ee
where
$$\Delta W^k=U^k_1 \sqrt{\Delta},\; \Delta Z^k={\textstyle\frac{1}{2}}\Delta^{3/2}(U^k_1+{\textstyle\frac{1}{\sqrt{3}}}U^k_2),$$
with $U^k_1,U^k_2$ being independent standard Gaussian random variables and where $\Delta$ is the time step. Moreover,
$$L^0=b_1\partial_1+b_2\partial_2+{\textstyle \frac{1}{2}}\left(\partial_1^2+\partial_2^2\right),\quad L^1=\partial_1,\, L_2=\partial_2.$$ The corresponding results are shown in Figure \ref{graphs} with $\Delta=10^{-4}$. Note that for $f=1_A$, $A$ a Borel set, equality (\ref{ergf}) entails
\be\label{ergl}
\lim_{T\rightarrow \infty}\frac{1}{T}\int_0^T 1_A(X_t)dt = \mu(A).
\ee
Hence indicator functions comprise a natural class of test functions to be used. A glance at Figure \ref{graphs} clearly displays a convergent behavior of expression (\ref{fT}) for different starting points and different positions of the indicator functions. One also sees that the limits decay rather fast and that the diffusion takes ever more time to stabilize the farther away we move from the origin. On the other hand, if convergence of (\ref{ergl}) is observed then the limit $\mu$ is a natural candidate for an invariant probability measure. Note that according to \cite[Theorem 2.2.5]{Lor} all transition probabilities $K_t(x,dy)$ of the diffusion $X$ are equivalent which by Doob's theorem, see \cite[Theorem 4.2.1]{DaP}, gives the uniqueness of the invariant measure $\mu$. This means that the numerical findings support the assertion on ergodicity of $X$.
\begin{figure}[!]
\begin{center}
\includegraphics[scale = 0.9]{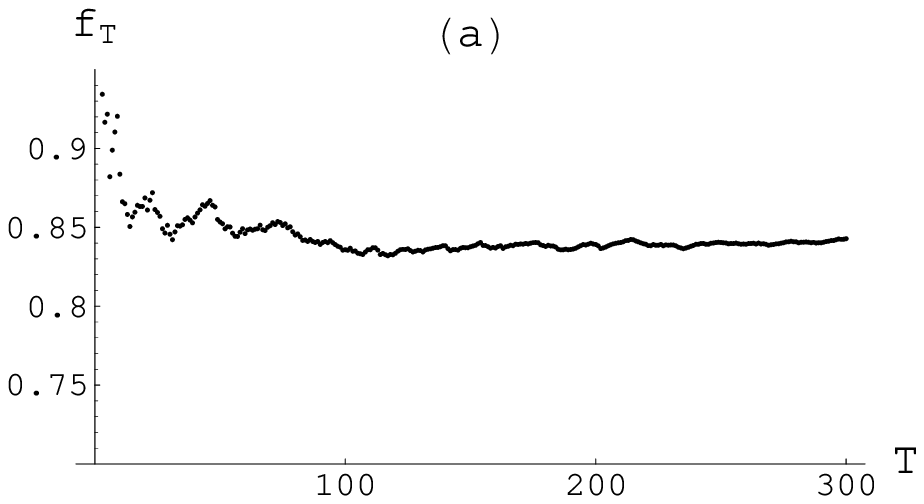}
\includegraphics[scale = 0.9]{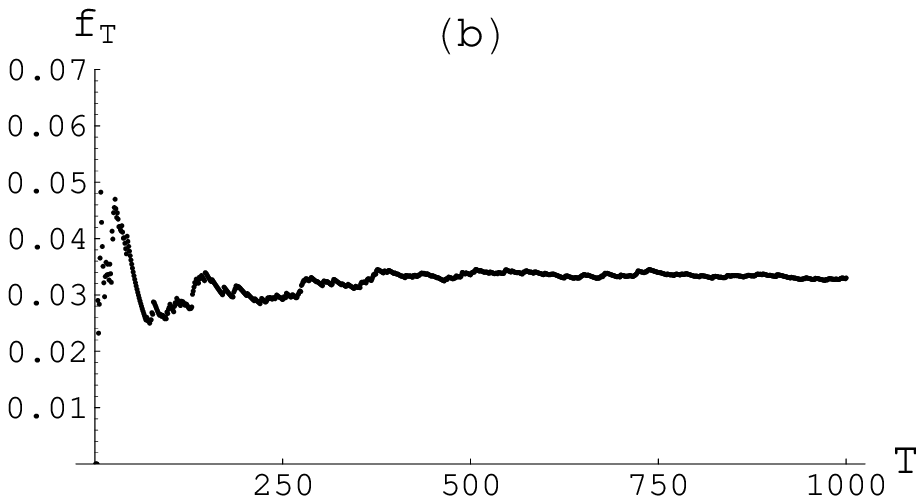}\\[2ex]
\includegraphics[scale = 0.9]{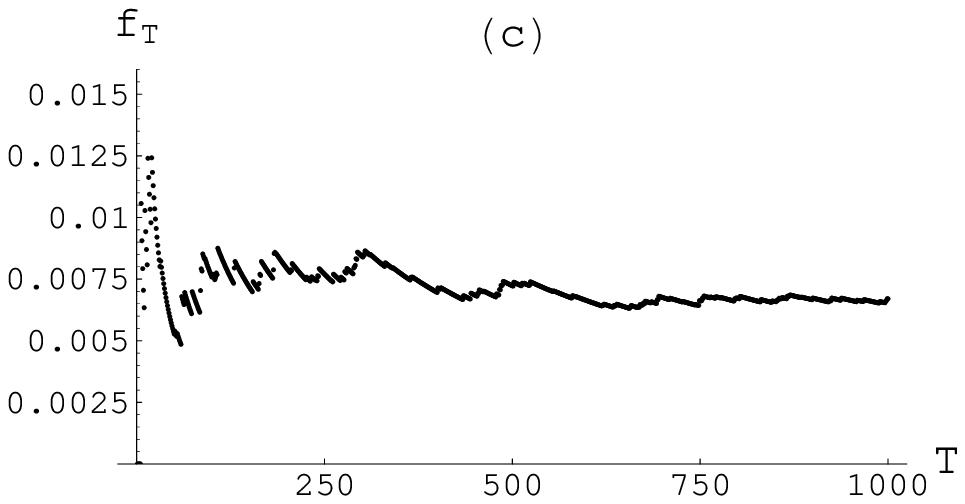}
\includegraphics[scale = 0.9]{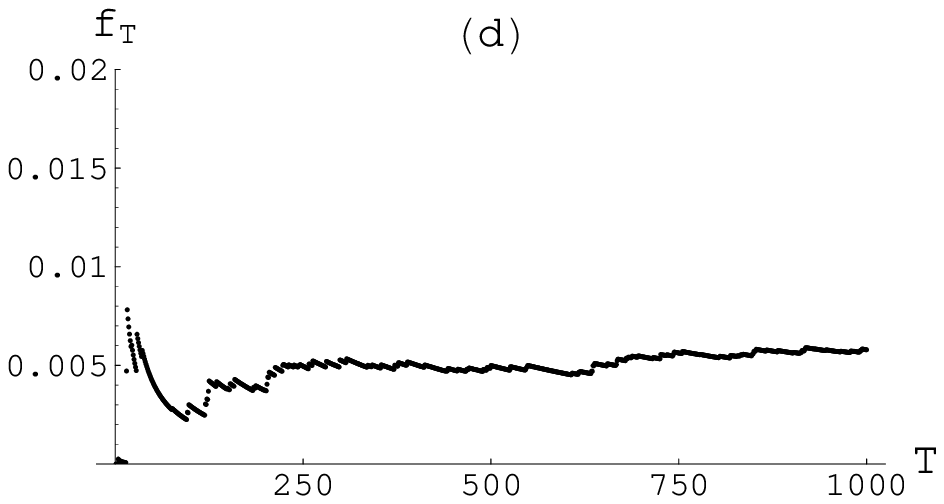}
\end{center}
\end{figure}
\section{Acknowledgments} The author would like to thank Simona Bernabei for her careful reading of the manuscript. He also wishes to express his gratitude for the financial support through the program ``Rientro dei cervelli'' of the italian M.I.U.R.

\end{document}